\documentclass[review]{elsarticle}

\usepackage{lineno,hyperref}
\usepackage{amssymb,enumerate}
\usepackage{amsmath}
\usepackage{graphics}
\newtheorem{thm}{Theorem}[section]
 
 \newtheorem{lem}[thm]{Lemma}

 \newtheorem{rem}[thm]{Remark}

\modulolinenumbers[5]

\journal{Journal of \LaTeX\ Templates}

\begin{document}

\begin{frontmatter}

\title{Blow-up and lifespan estimate to a nonlinear wave equation in Schwarzschild spacetime}


\author[mymainaddress,secondaryaddress]{Ning-An Lai\corref{mycorrespondingauthor}}
\cortext[mycorrespondingauthor]{Corresponding Author}
\ead{hyayue@gmail.com}

\address[mymainaddress]{Institute of Nonlinear Analysis and Department of Mathematics,\\ Lishui University, Lishui 323000, China}
\address[secondaryaddress]{School of Mathematical Sciences, Fudan University, Shanghai 200433, China}

\begin{abstract}

In \cite{Luk}, Luk proved global existence for semilinear wave equations in Kerr spacetime with small angular momentum($a\ll M$)
\[
\Box_{g_K}\phi=F(\partial \phi),
\]
when the quadratic nonlinear term satisfies the null condition. In this work, we will show that if the null condition does not hold, at most we can have almost global existence for semilinear wave equations with quadratic
nonlinear term in Schwarzschild spacetime, which is the special case of Kerr with $a=0$
\[
\Box_{g_S}\phi=(\partial_t\phi)^2,
\]
where $\Box_{g_S}$ denotes the D'Alembert operator associated with Schwarzschild metric. What is more, if the power of the nonlinear term is replaced with $p$ satisfying $3/2\le p< 2$, we still can show blow-up, no matter how small the initial data are. We do not have to assume that the support of the data should be far away from the event horizon.
\end{abstract}


\end{frontmatter}


\section{Introduction}

In this work we are going to study the nonlinear wave equation in Schwarzschild spacetime
\begin{equation}\label{equ}
\Box_{g_S} u=|u_t|^p,~~~in ~\mathbb{R}_t^+\times \Sigma,\\
\end{equation}
where $g_S$ denotes the Schwarzschild metric
\begin{equation}\label{2}
\begin{aligned}
g_S=F(r)dt^2-F(r)^{-1}dr^2-r^2d\omega^2
\end{aligned}
\end{equation}
with $F(r)=1-\frac{2M}{r}$, the constant $M>0$ has the meaning of mass of the black hole and $d\omega^2$ is the standard metric on the unit sphere $\mathbb{S}^2$.
The D'Alembert operator associated with the Schwarzschild metric $g$ has the form
\begin{equation}\label{3}
\begin{aligned}
\Box_{g_S}=\frac{1}{F(r)}\Big(\partial_t^2-\frac{F}{r^2}\partial_r(r^2F)\partial_r-\frac{F}{r^2}\Delta_{\mathbb{S}^2}\Big),
\end{aligned}
\end{equation}
where $\Delta_{\mathbb{S}^2}$ denotes the standard Laplace-Beltrami operator on $\mathbb{S}^2$. $\mathbb{R}_t^+\times \Sigma$ is called the
exterior of the black hole:
\begin{equation}\nonumber\\
\begin{aligned}
\mathbb{R}_t^+\times \Sigma=\mathbb{R}_t^+\times (2M, \infty)\times \mathbb{S}^2.\\
\end{aligned}
\end{equation}

Such Cauchy problem of semilinear wave equations with derivatives in the nonlinear term was first proposed by John \cite{John}, in which he studied the following quasilinear wave equation in $\mathbb{R}^3$
\begin{equation}\label{4}
\begin{aligned}
\partial_t^2u(t, x)-\Delta u(t, x)=\partial_t|au(t, x)+b\partial_tu(t, x)|^p
\end{aligned}
\end{equation}
and proved that the solutions blow up at a finite time when $p=2$ and the compact supported data satisfy some conditions.
His method is also applicable to the following Cauchy problem
\begin{equation}\label{5}
\begin{aligned}
\partial_t^2u(t, x)-\Delta u(t, x)=|u_t|^p
\end{aligned}
\end{equation}
with $1<p\leq 2$. Assuming the data is radial symmetric and $p>2$, Sideris \cite{Sideris} showed that the Cauchy problem \eqref{5} in $\mathbb{R}^3$ admits global small amplitude solutions. Klainerman \cite{Klainerman1} and \cite{Klainerman2} proved global existence of small data solutions to nonlinear wave equations of a more general form
\begin{equation}\label{6}
\left \{
\begin{aligned}
&\partial_t^2u(t, x)-\Delta u(t, x)=F(\partial_tu, \nabla_xu),\\
&|D^{\alpha}F(w)|\leq A|w|^{p-\alpha},~~p>\alpha,~|w|\leq 1
\end{aligned} \right.
\end{equation}
with $p>\frac{n+1}{n-1}$. Hence it is easy to see that problem \eqref{5} in $\mathbb{R}^3$ admits a critical power $p=2$.
In $\mathbb{R}^2$, Schaeffer \cite{Schaeffer} showed blow-up result for problem \eqref{5} with $p=3$ and conjectured that it is the critical exponent(see also John \cite{John2}). Later, Agemi \cite{Agemi} proved that problem \eqref{5} in $\mathbb{R}^2$ with $1<p\leq 3$ has no global solutions. The case $n=1$ for problem \eqref{5} is essentially due to Masuda \cite{Masuda}, in which he established blow-up result in $\mathbb{R}^n(n=1, 2, 3)$ with power $p=2$. Hidano and
Tsutaya \cite{HT} and Tzvetkov \cite{Tzv} proved global existence for $p>p_c(n)$ and $n=2, 3$ independently.
Zhou \cite{Zhou1} demonstrated a simple proof of blow-up result to problem \eqref{5} in $\mathbb{R}^n$ with $1<p\leq \frac{n+1}{n-1}$. Recently, Hidano, Wang and Yokoyama \cite{HWY} established global existence for $p>p_c(n)$ and $n\ge 4$, by assuming the data are radially symmetrical. We also refer the reader to a recent paper \cite{LT} for semilinear wave equations related to Glassey's conjecture with scattering damping.

Recently, Zhou and Han \cite{Zhou2} established blow-up result to the initial boundary value problem \eqref{5} in exterior domain with $1<p\leq \frac{n+1}{n-1}(n\geq 1)$, by introducing a test function. Also they obtained the upper bound of lifespan estimate. Such test function method can be dated back to Yordanov and Zhang \cite{Yordanov}, in which they studied the Strauss conjecture with critical power in $\mathbb{R}^n(n\geq 4)$.

In the Schwarzschild spacetime, Catania and Georgiev \cite{Catania} considered the Cauchy problem of semilinear wave equation
\[
 \Box_{g_S} u=|u|^p.
\]
 In the Regge-Wheeler coordinates, by choosing special initial data they proved the blow-up results in two cases: (I) $1<p<1+\sqrt 2$ and small data supported far away from the event horizon; (II) $2<p<1+\sqrt 2$ and large data supported near the event horizon. On the contrary, Lindblad et al \cite{LMS} established global existence for $p>1+\sqrt2$, which also holds for Kerr black hole with small angular momentum if the small data have compact support.
 One can also find the global existence result for $p>3$ in \cite{BS} and $p=5$ in \cite{MMT}.

For a black hole metric $g$, one of the key ingredients to prove asymptotic stability is to understand the nonlinear toy problem
\begin{equation}\label{toy}
\Box_g \phi=(\partial \phi)^2,
\end{equation}
especially when the quadratic nonlinear term satisfies some special structure, i.e., null condition, due to the reason that null condition has served as a good model problem for the study of stability of the flat Euclidean spacetime.
Luk \cite{Luk} studied the semilinear wave equations on a Kerr spacetime with small angular momentum($a\ll M$)
\[
\Box_{g_K}\phi=F(\partial \phi),
\]
where $g_K$ denotes the Kerr metric, which is  parametrized by two parameters $(M, a)$, representing the mass and angular momentum of a black hole, respectively. And the nonlinear term $F$ is at least quadratic. Assuming the nonlinear term $F$ satisfies the null condition, then he can prove global existence for any initial data that are sufficiently small. In this work we are devoted to studying the Cauchy problem of semiliner wave equation with derivatives in Schwarzschild spacetime, a special example of Kerr black hole with $a=0$, as stated in \eqref{equ}.
We will show that if the quadratic nonlinear term does not satisfy the null condition, the the solution will blow up in a finite time, no matter how small the initial data are. Also, the lifespan estimate of exponential type will be established. What is more, if $\frac 32\le p<2$, we then show that the lifespan is of polynomial type.
The first key step of our proof is to rewrite the equation in the Regge-Wheeler coordinates, and then choose the test function which was introduced in \cite{Catania}. The second one is to construct an auxiliary functional which satisfies an ordinary differential inequality of Riccati type, following the idea of \cite{Zhou2}. We mention that the approach in  \cite{Catania} is based on applying variants of the classical Kato's lemma to an auxiliary quantity, constructed from the corresponding solutions and test function.

\begin{rem}
We do not have to assume that the support of the initial data should be far away from the event horizon.
\end{rem}

\begin{rem}
It is interesting to verify that $p=2$ is the critical exponent for the Cauchy problem \eqref{equ}, thus, to show global existence for $p>2$.
\end{rem}

\begin{rem}
Our method does not work for the smaller power $1<p<\frac 32$, due to the singularity when approaching to the event horizon.
\end{rem}

As stated above, our main goal is to show blow-up to the Cauchy problem of semilinear wave equations in Schwarzschild spacetime for $\frac 32\leq p\leq 2$. We first introduce the Regge-Wheeler coordinate
\begin{equation}\label{7a}
\begin{aligned}
s(r)=r+2M\ln(r-2M),
\end{aligned}
\end{equation}
then the equation in \eqref{equ} can be written as
\begin{equation}\label{7}
\begin{aligned}
\partial_t^2u-\partial_s^2u-\frac{2F}{r}\partial_su-\frac{F}{r^2}\Delta_{\mathbb{S}^2}u=F|u_t|^p,
\end{aligned}
\end{equation}
where
\begin{equation}\nonumber\\
\begin{aligned}
F=F(s)=1-\frac{2M}{r(s)}
\end{aligned}
\end{equation}
and $r=r(s)$ is the inverse function of \eqref{7a}.

Without loss of generality, in what follows we consider the radial solutions, thus solutions of the form: $u=u(t, s)$. Hence equation \eqref{7} is simplified to
\begin{equation}\label{8}
\begin{aligned}
\partial_t^2u-\partial_s^2u-\frac{2F}{r}\partial_su=F|u_t|^p.
\end{aligned}
\end{equation}
Set $v(t, s)=r(s)u(t, s)$, we then can rewrite equation \eqref{8} in the following form
\begin{equation}\label{9}
\begin{aligned}
\partial_t^2v(t, s)-\partial_s^2v(t, s)+\frac{2MF}{r^3}v(t, s)=Fr^{1-p}|v_t|^p.
\end{aligned}
\end{equation}

Noting that $s$ varies from $-\infty$ to $+\infty$ as $r$ varies from $2M$ to $+\infty$, then if we restrict ourselves to symmetrically radial solutions, we may focus on the following Cauchy problem
\begin{equation}\label{10}
\left \{
\begin{aligned}
&\partial_t^2v(t, s)-\partial_s^2v(t, s)+W(s)v(t, s)=h(s)|v_t|^p,~(t, s)\in [0, \infty)\times \mathbb{R},\\
&v(0, s)=\varepsilon f(s),~~ v_t(0, s)=\varepsilon g(s),~~s\in \mathbb{R},\\
\end{aligned} \right.
\end{equation}
where $h(s)=F(s)r(s)^{1-p}$ and $W(s)=\frac{2MF(s)}{r(s)^3}$ with $F(s)=1-\frac{2M}{r(s)}$. The initial data $f(s), g(s)$ are nonnegative satisfying that $g(s)$ does not vanish and
\[
supp~f, g\subset \{s:\big||s|\le R\}
\]
with some constant $R>0$.

Our main result reads
\begin{thm}\label{thm1}
Let $\frac32\leq p\leq 2$, then the solutions of problem \eqref{10} blow up in a finite time $T(\varepsilon)$. Furthermore, we obtain the upper bound of lifespan estimate as follows.\\
(I) If $\frac32\le p<2$, then there exists a positive constant $C_1$ which is independent of $\varepsilon$ such that
 \begin{equation}\label{13}
\begin{aligned}
T(\varepsilon)\leq C_1\varepsilon^{-\frac{p-1}{2-p}}.
\end{aligned}
\end{equation}
(II) If $p=2$, then there exists a positive constant $C_2$ which is independent of $\varepsilon$ such that
 \begin{equation}\label{14}
\begin{aligned}
T(\varepsilon)\leq \exp\left(C_2\varepsilon^{-1}\right).
\end{aligned}
\end{equation}
\end{thm}

We arrange the rest of the paper as follows. In section 2 we give some preliminary lemmas. In section 3 we demonstrate the proof of Theorem \ref{thm1}.
\section{Preliminaries}
According to Birkhoff's theorem, the Schwarzschild metric is the most general spherically symmetric, vacuum solution of the Einstein's field equations. We refer the reader to \cite{Kong} for a general introduction of the Einstein's field equations and its exact solutions. The sphere $\{r=2M\}$, referred to as the event horizon, is merely a coordinate singularity, while the origin $\{r=0\}$ denotes a true curvature singularity.

Set
\begin{equation}\label{14a}
\begin{aligned}
h(s)=F(s)r(s)^{1-p},
\end{aligned}
\end{equation}
then from the definition of \eqref{7a} we have the following asymptotic behavior.
\begin{lem}\label{lem3}
Let $h(s)$ be as in \eqref{14a}, then it holds
\begin{equation}\label{14b}
h(s)\sim\left \{
\begin{aligned}
&s^{1-p},~~s\geq 4M+e,\\
&e^{\frac{s}{2M}},~~-\infty<s\leq 4M+e,
\end{aligned} \right.
\end{equation}
where the notation $X\sim Y$ means $C^{-1}Y\leq X\leq CY$ with some positive constant $C$, and $e$ denotes the Euler's number.
\end{lem}
\textit{\emph Proof.} By definition \eqref{7a}, it is easy to get
\[
s\ge 4M+e \Leftrightarrow r\ge 2M+e,
\]
which yields
\begin{equation}\label{Fs}
\begin{aligned}
\frac{e}{4M+e}\le 1-\frac{2M}{r}\le 1.
\end{aligned}
\end{equation}
Also, if $r\ge 2M+e$, then there exists a positive constant such that
\begin{equation}\label{rs}
\begin{aligned}
r\le s\le Cr.
\end{aligned}
\end{equation}
And hence the first part of \eqref{14b} follows from \eqref{Fs} and \eqref{rs}.

On the other hand, one has
\[
s\le 4M+e\Leftrightarrow 2M\le r\le 2M+e,
\]
which means that
\begin{equation}\label{hs}
\begin{aligned}
h(s)=(r-2M)r^{-p}\sim C(r-2M).
\end{aligned}
\end{equation}
By \eqref{7a}, we know that for some positive constant $C$
\[
2M\ln(r-2M)=s-r\sim s-C,
\]
which yields
\begin{equation}\label{ln}
\begin{aligned}
r-2M\sim e^{\frac{s-C}{2M}}=e^{-\frac{C}{2M}}e^{\frac{s}{2M}}.
\end{aligned}
\end{equation}
Therefor, the second part of \eqref{14b} comes from \eqref{hs} and \eqref{ln}.

\begin{lem}(\cite{Catania}, Lemma 5.3)\label{lem1}
Let $W(s)=\frac{2MF(s)}{r(s)^3}$ with $F(s)=1-\frac{2M}{r(s)}$, then given any $A>0$, the equation
\begin{equation}\label{16}
\begin{aligned}
\big(-\partial_s^2+W(s)+A^2\big)\varphi(s)=0,~~s\in \mathbb{R}
\end{aligned}
\end{equation}
admits a positive solution $\varphi\in C^2(\mathbb{R})$ such that
\begin{equation}\label{17}
\begin{aligned}
\varphi(s)\thicksim e^{As},~~|s|\rightarrow \infty.
\end{aligned}
\end{equation}
\end{lem}

\begin{rem}
Noting that if the term $W(s)$ disappears in \eqref{16}, then $e^{As}$ is an exact solution.
\end{rem}

\section{Proof of Theorem \ref{thm1}}

We are now in a position to show the proof of Theorem \ref{thm1}. First we introduce a lemma which will be used in our proof. In the following $C$ denotes a positive constant which may change from line to line.
\begin{lem}\label{lem2}
Given any $\alpha\geq 0, \beta>0$ and $L>0$, there exists a positive constant $C$ such that
\begin{equation}\label{16a}
\begin{aligned}
\int_{0<s\leq t+L}(1+s)^{\alpha}e^{-\beta(t-s)}ds\leq C(t+L)^{\alpha}.
\end{aligned}
\end{equation}
\end{lem}
\textit{\emph Proof.} We split the integral in \eqref{16a} into two parts
\begin{equation}\label{16b}
\begin{aligned}
&\int_{0<s\leq t+L}(1+s)^{\alpha}e^{-\beta(t-s)}ds\\
=&\int_{0}^{\frac{t+L}{2}}(1+s)^{\alpha}e^{-\beta(t-s)}ds+\int_{\frac{t+L}{2}}^{t+L}
(1+s)^{\alpha}e^{-\beta(t-s)}ds\\
\leq &Ce^{-\beta\cdot \frac{t-L}{2}}\int_{0}^{\frac{t+L}{2}}(1+s)^{\alpha}ds+C(t+L)^{\alpha}\int_{\frac{t+L}{2}}^{t+L}
e^{-\beta(t-s)}ds\\
\leq &C(t+L)^{\alpha},\\
\end{aligned}
\end{equation}
which is the desired inequality \eqref{16a}.

According to Lemma \ref{lem1}, the following equation
\begin{equation}\label{18}
\begin{aligned}
\big(-\partial_s^2+W(s)+\frac{1}{4M^2}\big)\phi(s)=0,~~s\in \mathbb{R},
\end{aligned}
\end{equation}
admits a positive solution $\phi(s)$ satisfying
\begin{equation}\label{19}
\begin{aligned}
\phi(s)\thicksim e^{\frac{s}{2M}},~~~|s|\rightarrow \infty.
\end{aligned}
\end{equation}

Set
\[
 \psi(t, s)=e^{-\frac{t}{2M}}\phi(s),
 \]
then it is easy to see that
\begin{equation}\label{20}
\begin{aligned}
&\psi_t(t, s)=-\frac{1}{2M}\psi(t, s),\\
&\partial_s^2\psi(t, s)-W(s)\psi(t, s)=\frac{1}{4M^2}\psi(t, s).
\end{aligned}
\end{equation}
Multiplying the both sides of equation in \eqref{10} with $\psi(t, s)$ and integrating by parts yields
\begin{equation}\label{20a}
\begin{aligned}
&\int_{\mathbb{R}}\psi(t, s)\big(\partial_t^2v-\partial_s^2v+W(s)v\big)ds\\
=&\int_{\mathbb{R}}\psi\partial_t^2v-\big(\partial_s^2\psi-W(s)\psi\big)vds\\
=&\int_{\mathbb{R}}\Big(\psi\partial_t^2v-\frac{1}{4M^2}\psi v\Big)ds\\
=&\int_{\mathbb{R}}h(s)|v_t|^p\psi ds,\\
\end{aligned}
\end{equation}
where the second equality in \eqref{20} has been used.

By direct calculation we get
\begin{equation}\label{21}
\begin{aligned}
\frac{d}{dt}\int_{\mathbb{R}}\psi v_tds&=\int_{\mathbb{R}}\psi v_{tt}ds+\int_{\mathbb{R}}\psi_tv_tds\\
&=\int_{\mathbb{R}}\big(\psi v_{tt}-\frac{1}{2M}\psi v_t\big)ds
\end{aligned}
\end{equation}
and
\begin{equation}\label{22}
\begin{aligned}
\frac{d}{dt}\int_{\mathbb{R}}\psi vds&=\int_{\mathbb{R}}\psi v_{t}ds+\int_{\mathbb{R}}\psi_tvds\\
&=\int_{\mathbb{R}}\big(\psi v_{t}-\frac{1}{2M}\psi v\big)ds.
\end{aligned}
\end{equation}
Multiplying \eqref{22} with $\frac{1}{2M}$, adding it to \eqref{21} and then combining with \eqref{20} we come to
\begin{equation}\label{23}
\begin{aligned}
&\frac{d}{dt}\int_{\mathbb{R}}\psi v_tds+\frac{1}{2M}\frac{d}{dt}\int_{\mathbb{R}}\psi vds\\
=&\int_{\mathbb{R}}\Big(\psi v_{tt}-\frac{1}{4M^2}\psi v\Big)ds\\
=&\int_{\mathbb{R}}h(s)|v_t|^p\psi ds.\\
\end{aligned}
\end{equation}
It follows that if we integrating \eqref{23} over $[0, t]$
\begin{equation}\label{24}
\begin{aligned}
&\int_{\mathbb{R}}\big(\psi v_t+\frac{1}{2M}\psi v\big)ds\\
=&\int_{\mathbb{R}}\phi(s)g(s)ds+\frac{1}{2M}\int_{\mathbb{R}}\phi(s)f(s)ds+\int_0^t\int_{\mathbb{R}}h(s)|v_{\tau}|^p\psi dsd\tau\\
\geq &\int_{\mathbb{R}}\phi(s)g(s)ds+\int_0^t\int_{\mathbb{R}}h(s)|v_{\tau}|^p\psi dsd\tau,\\
\end{aligned}
\end{equation}
where we have used the fact that the initial data are nonnegative.

Multiplying \eqref{24} with $\frac{1}{2M}$ and then adding it to \eqref{20}, one gets
\begin{equation}\label{25}
\begin{aligned}
&\int_{\mathbb{R}}\big(\psi v_{tt}+\frac{1}{2M}\psi v_t\big)ds\\
\geq &\int_{\mathbb{R}}h(s)|v_{t}|^p\psi ds+\frac{1}{2M}\int_{\mathbb{R}}\phi(s)g(s)ds+\frac{1}{2M}\int_0^t\int_{\mathbb{R}}h(s)|v_{\tau}|^p\psi dsd\tau.\\
\end{aligned}
\end{equation}
Hence we have
\begin{equation}\label{26}
\begin{aligned}
&\frac{d}{dt}\int_{\mathbb{R}}\psi v_tds+\frac{1}{M}\int_{\mathbb{R}}\psi v_tds\\
=&\int_{\mathbb{R}}\psi v_{tt}ds-\frac{1}{2M}\int_{\mathbb{R}}\psi v_tds+\frac{1}{M}\int_{\mathbb{R}}\psi v_tds\\
=&\int_{\mathbb{R}}\big(\psi v_{tt}+\frac{1}{2M}\psi v_t\big)ds\\
\geq &\int_{\mathbb{R}}h(s)|v_{t}|^p\psi ds+\frac{1}{2M}\int_{\mathbb{R}}\phi(s)g(s)ds+\frac{1}{2M}\int_0^t\int_{\mathbb{R}}h(s)|v_{\tau}|^p\psi dsd\tau.\\
\end{aligned}
\end{equation}

Set
\begin{equation}\label{27}
\begin{aligned}
G(t)=\int_{\mathbb{R}}\psi v_tds-\frac{1}{2}\int_0^t\int_{\mathbb{R}}h(s)|v_{\tau}|^p\psi dsd\tau-\frac{1}{2}\int_{\mathbb{R}}\phi(s)g(s)ds.
\end{aligned}
\end{equation}
Then we have
\begin{equation}\label{28}
\begin{aligned}
G(0)=\frac12\int_{\mathbb{R}}\phi(s)g(s)ds\geq 0,
\end{aligned}
\end{equation}
and
\begin{equation}\label{29}
\begin{aligned}
\frac{dG(t)}{dt}=\int_{\mathbb{R}}\Big(\psi v_{tt}-\frac{1}{2M}\psi v_t\Big)ds-\frac{1}{2}\int_{\mathbb{R}}h(s)|v_{t}|^p\psi ds.
\end{aligned}
\end{equation}
By combining \eqref{26}, \eqref{27} and \eqref{29} we arrive at
\begin{equation}\label{30}
\begin{aligned}
&\frac{dG(t)}{dt}+\frac{1}{M}G(t)\\
=&\int_{\mathbb{R}}\Big(\psi v_{tt}+\frac{1}{2M}\psi v_t\Big)ds-\frac{1}{2}\int_{\mathbb{R}}h(s)|v_{t}|^p\psi ds\\
&-\frac{1}{2M}\int_0^t\int_{\mathbb{R}}h(s)|v_{\tau}|^p\psi dsd\tau-\frac{1}{2M}\int_{\mathbb{R}}\phi(s)g(s)ds\\
\geq &\int_{\mathbb{R}}h(s)|v_{t}|^p\psi ds+\frac{1}{2M}\int_{\mathbb{R}}\phi(s)g(s)ds+\frac{1}{2M}\int_0^t\int_{\mathbb{R}}h(s)|v_{\tau}|^p\psi dsd\tau\\
&-\frac{1}{2}\int_{\mathbb{R}}h(s)|v_{t}|^p\psi ds-\frac{1}{2M}\int_0^t\int_{\mathbb{R}}h(s)|v_{\tau}|^p\psi dsd\tau\\
&-\frac{1}{2M}\int_{\mathbb{R}}\phi(s)g(s)ds\\
=&\frac{1}{2}\int_{\mathbb{R}}h(s)|v_{t}|^p\psi ds,\\
\end{aligned}
\end{equation}
which implies
\begin{equation}\label{31}
\begin{aligned}
\frac{d}{dt}\Big(e^{\frac{t}{M}}G(t)\Big)\geq 0,
\end{aligned}
\end{equation}
and since $G(0)\geq 0$ this in turn gives
\begin{equation}\label{32}
\begin{aligned}
G(t)=\int_{\mathbb{R}}\psi v_tds-\frac{1}{2}\int_0^t\int_{\mathbb{R}}h(s)|v_{\tau}|^p\psi dsd\tau-\frac{1}{2}\int_{\mathbb{R}}\phi(s)g(s)ds\geq 0.
\end{aligned}
\end{equation}

Set
\begin{equation}\label{33}
\begin{aligned}
F(t)=\frac{1}{2}\int_0^t\int_{\mathbb{R}}h(s)|v_{\tau}|^p\psi dsd\tau+\frac{1}{2}\int_{\mathbb{R}}\phi(s)g(s)ds\geq 0,
\end{aligned}
\end{equation}
then we conclude from \eqref{32} that
\begin{equation}\label{34}
\begin{aligned}
F(t)\leq \int_{\mathbb{R}}\psi v_tds.
\end{aligned}
\end{equation}
By finite speed of propagation property, H\"{o}lder inequality and \eqref{19} we have
\begin{equation}\label{35}
\begin{aligned}
\Big|\int_{\mathbb{R}}\psi v_tds\Big|^p&=\Big|\int_{|s|\leq t+R}\psi v_tds\Big|^p\\
&\leq \int_{\mathbb{R}}h|v_t|^p\psi ds\Big(\int_{|s|\leq t+R}h^{-\frac{1}{p-1}}\psi ds\Big)^{p-1}\\
&\leq C\int_{\mathbb{R}}h|v_t|^p\psi ds\Big(\int_{|s|\leq t+R}h^{-\frac{1}{p-1}}e^{\frac{s-t}{2M}} ds\Big)^{p-1}.\\
\end{aligned}
\end{equation}
Let
\[
I\triangleq \int_{|s|\leq t+R}h^{-\frac{1}{p-1}}e^{\frac{s-t}{2M}} ds,
\]
noting the asymptotic behavior \eqref{14b}, we then may estimate $I$ as
\begin{equation}\label{I}
\begin{aligned}
I&=\int_{-t-R}^{4M+e}h^{-\frac{1}{p-1}}e^{\frac{s-t}{2M}} ds+\int_{4M+e}^{t+R}h^{-\frac{1}{p-1}}e^{\frac{s-t}{2M}} ds\\
&\le C\int_{-t-R}^{4M+e}e^{-\frac{s}{2M}\frac{1}{p-1}}e^{\frac{s-t}{2M}} ds+C\int_{4M+e}^{t+R}(1+s)e^{\frac{s-t}{2M}} ds\\
&\le II+C(t+R),
\end{aligned}
\end{equation}
where we have used Lemma \ref{lem2} with $\alpha=1$ and $\beta=\frac{1}{2M}$ in the last inequality, and
\[
II=C\int_{-t-R}^{4M+e}e^{-\frac{s}{2M}\frac{1}{p-1}}e^{\frac{s-t}{2M}} ds.
\]

Next we estiamte $II$ in two cases: $\frac32\leq p<2$ and $p=2$.
\subsection{$\frac32\leq p<2$}
In this case, by direct computation we hvae
\begin{equation}\label{35a}
\begin{aligned}
II&=Ce^{-\frac{t}{2M}}\int_{-t-R}^{4M+e}e^{\frac{s}{2M}\frac{p-2}{p-1}}e^{\frac{s-t}{2M}} ds\\
&\le Ce^{-\frac{t}{2M}\frac{2p-3}{p-1}}\\
&\le C(t+R).
\end{aligned}
\end{equation}
\subsection{$p=2$}
In this case, it is easy to get
\begin{equation}\label{35aa}
\begin{aligned}
II&=Ce^{-\frac{t}{2M}}\int_{-t-R}^{4M+e}ds\\
&\le C(t+R).
\end{aligned}
\end{equation}
We then conclude from \eqref{I}, \eqref{35a} and \eqref{35aa} that for $\frac32\le p\le 2$
\begin{equation}\label{I1}
\begin{aligned}
I<C(t+R).
\end{aligned}
\end{equation}

Then by combining definition \eqref{33}, inequality \eqref{34}, \eqref{35} and \eqref{I1} it follows that for $\frac32\le p\leq 2$
\begin{equation}\label{36}
\begin{aligned}
F'(t)&=\frac12\int_{\mathbb{R}}f|v_t|^p\psi ds\\
&\geq \frac{C\Big|\int_{\mathbb{R}}\psi v_tds\Big|^p}{(t+R)^{p-1}}\\
&\geq \frac{C|F(t)|^p}{(t+R)^{p-1}}.\\
\end{aligned}
\end{equation}
Hence by the property of Riccati equation, $F(t)$ must blow up at a finite time if $p-1\leq 1$, thus: $p\leq 2$. We finish the blow-up part of Theorem \ref{thm1}.

Next we will establish the upper bound of the lifespan estimates. It is easy to see that $F(t)$ satisfies
the following problem
\begin{equation}\label{37}
\left \{
\begin{aligned}
&F'(t)\geq \frac{C|F(t)|^p}{(t+R)^{p-1}},\\
&F(0)=N\varepsilon,
\end{aligned}\right.
\end{equation}
with
\begin{equation}\nonumber\\
\begin{aligned}
N=\frac12\int_{\mathbb{R}}\phi(s)g(s)ds>0.
\end{aligned}
\end{equation}
If we consider the Riccati equation
\begin{equation}\label{38}
\left \{
\begin{aligned}
&H'(t)= \frac{C|H(t)|^p}{(t+R)^{p-1}},\\
&H(0)=N\varepsilon,
\end{aligned} \right.
\end{equation}
then $F(t)$ blows up before $H(t)$, which means that we can estimate the upper bound of lifespan of $F(t)$ through that of $H(t)$.

In the case $\frac32\le p<2$, we solve problem \eqref{38} and get
\begin{equation}\label{39}
\begin{aligned}
H(t)=\Big((N\varepsilon)^{1-p}+\widetilde{C}R^{2-p}-\widetilde{C}(t+R)^{2-p}\Big)^{-\frac{1}{p-1}}
\end{aligned}
\end{equation}
with $\widetilde{C}=\frac{C(p-1)}{2-p}$. Therefore we can estimate the lifespan of $F(t)$ as
\begin{equation}\label{40}
\begin{aligned}
T(\varepsilon)\leq C_1\varepsilon^{-\frac{p-1}{2-p}}
\end{aligned}
\end{equation}
with
\begin{equation}\nonumber\\
\begin{aligned}
C_1=\Big(\frac{2-p}{C(p-1)}N^{1-p}\Big)^{\frac{1}{2-p}}>0,
\end{aligned}
\end{equation}
which is independent of $\varepsilon$.

In the case $p=2$, we solve problem \eqref{38} and get
\begin{equation}\label{41}
\begin{aligned}
H(t)=\Big((N\varepsilon)^{-1}-C\ln{\frac{t+R}{R}}\Big)^{-1},
\end{aligned}
\end{equation}
from which we derive the lifespan from above of $F(t)$ as
\begin{equation}\label{40a}
\begin{aligned}
T(\varepsilon)\leq Ce^{C_2\varepsilon^{-1}}
\end{aligned}
\end{equation}
with
\begin{equation}\nonumber\\
\begin{aligned}
C_2=\frac{1}{CN}>0,
\end{aligned}
\end{equation}
which is independent of $\varepsilon$. And hence we finish the proof of Theorem \ref{thm1}.

\section{Acknowledgement}

The author is supported by Natural Science Foundation of Zhejiang Province(LY18A010008), NSFC(11501273, 11726612), Postdoctoral Research Foundation of China(2017M620128, 2018T110332), the Scientific Research Foundation of the First-Class Discipline of Zhejiang Province
(B)(201601).

\section*{References}

\bibliography{mybibfile}

\end{document}